\documentclass[11pt]{amsart}
\usepackage{amssymb}

\usepackage{amsmath}
\usepackage[active]{srcltx}
\usepackage{t1enc}
\usepackage[latin2]{inputenc}
\usepackage{verbatim}
\usepackage{amsmath,amsfonts,amssymb,amsthm}
\usepackage[mathcal]{eucal}
\usepackage{enumerate}
\usepackage[centertags]{amsmath}
\usepackage{graphics}

\newtheorem{theorem}{Theorem}

\newtheorem{lemma}{Lemma}

\newtheorem{corollary}{Corollary}

\begin{document}
\author{George Tephnadze}
\title[Fej\'er means ]{On The maximal operators of Vilenkin-Fej\'er means}
\address{G. Tephnadze, department of Mathematics, Faculty of Exact and Natural
Sciences, Tbilisi State University, Chavchavadze str. 1, Tbilisi 0128,
Georgia}
\email{giorgitephnadze@gmail.com}
\date{}
\maketitle

\begin{abstract}
The main aim of this paper is to prove that the maximal operator $\overset{%
\sim }{\sigma }^{*}f:=\underset{n\in P}{\sup }\frac{\left| \sigma
_{n}f\right| }{\log ^{2}\left( n+1\right) }$ is bounded from the Hardy space
$H_{1/2}$ to the space $L_{1/2}$, where $\sigma _{n}f$ is Fejér means of
bounded Vilenkin-Fourier series.
\date{}
\end{abstract}



\section{ Introduction}

In one-dimensional case the weak type inequality
\begin{equation*}
\mu \left( \sigma ^{*}f>\lambda \right) \leq \frac{c}{\lambda }\left\|
f\right\| _{1},\text{ \qquad }\left( \lambda >0\right)
\end{equation*}
can be found in Zygmund \cite{Zy} for the trigonometric series, in Schipp
\cite{Sc} for Walsh series and in Pál, Simon \cite{PS} for bounded Vilenkin
series. Again in one-dimensional, Fujji \cite{Fu} and Simon \cite{Si2}
verified that $\sigma ^{*}$ is bounded from $H_{1}$ to $L_{1}$. Weisz \cite
{We2} generalized this result and proved the boundedness of $\sigma ^{*}$
from the martingale space $H_{p}$ to the space $L_{p}$ for $p>1/2$. Simon
\cite{Si1} gave a counterexample, which shows that boundedness does not hold
for $0<p<1/2.$ The counterexample for $p=1/2$ due to Goginava \cite{GoPubl},
(see also \cite{BGG2}). In the endpoint case $p=1/2$ two positive results
was showed. Weisz \cite{we4} proved that $\sigma ^{*}$ is bounded from the
Hardy space $H_{1/2}$ to the space weak-$L_{1/2}$. For Walsh-Paley system in
2008 Goginava \cite{GoSzeged} proved that the maximal operator $\widetilde{%
\sigma }$ $^{*\,}$defined by
\begin{equation*}
\widetilde{\sigma }^{*}f:=\sup_{n\in P}\frac{\left| \sigma _{n}f\right| }{%
\log ^{2}\left( n+1\right) }
\end{equation*}
is bounded from the Hardy space $H_{1/2}$ to the space $L_{1/2}$. He also
proved that for any nondecreasing function $\varphi :P_{+}\rightarrow [1,$ $%
\infty )$ satisfying the condition
\begin{equation}
\overline{\lim_{n\rightarrow \infty }}\frac{\log ^{2}\left( n+1\right) }{%
\varphi \left( n\right) }=+\infty  \label{cond}
\end{equation}
the maximal operator
\begin{equation*}
\sup_{n\in P}\frac{\left| \sigma _{n}f\right| }{\varphi \left( n\right) }
\end{equation*}
is not bounded from the Hardy space $H_{1/2}$ to the space $L_{1/2}.$

For Walsh-Kaczmarz system analogical theorem is proved in \cite{GNCz}.

The main aim of this paper is to prove that the maximal operator $\widetilde{%
\sigma }^{*}f$\thinspace with respect to Vilenkin system is bounded from the
Hardy space $H_{1/2}$ to the space $L_{1/2}$ (see Theorem 1).

We also prove that under the condition (\ref{cond}) the maximal operator
\begin{equation*}
\sup_{n\in P}\frac{\left| \sigma _{n}f\right| }{\varphi \left( n\right) }
\end{equation*}
is not bounded from the Hardy space $H_{1/2}$ to the space $L_{1/2}.$
Actually, we prove stronger result (see Theorem 2) than the unboundedness of
the maximal operator $\widetilde{\sigma }^{*}f$ from the Hardy space $%
H_{1/2} $ to the spaces $L_{1/2}.$ In particular, we prove that

\begin{equation*}
\underset{n\in P}{\sup }\left\| \frac{\sigma _{n}f}{\varphi \left( n\right) }%
\right\| _{L_{1/2}}=\infty .
\end{equation*}

\section{Definitions and Notation}

Let $P_{+}$ denote the set of the positive integers , $P:=P_{+}\cup \{0\}.$

Let $m:=(m_{0,}m_{1....})$ denote a sequence of the positive integers not
less than 2.

Denote by
\begin{equation*}
Z_{m_{k}}:=\{0,1,...m_{k}-1\}
\end{equation*}
the additive group of integers modulo $m_{k}.$

Define the group $G_{m}$ as the complete direct product of the group $%
Z_{m_{j}}$ with the product of the discrete topologies of $Z_{m_{j}}$ $^{,}$%
s.

The direct product $\mu $ of the measures
\begin{equation*}
\mu _{k}\left( \{j\}\right) :=1/m_{k},\text{ \qquad }(j\in Z_{m_{k}}),
\end{equation*}
is the Haar measure on $G_{m_{\text{ }}},$ with $\mu \left( G_{m}\right) =1.$

If $\sup\limits_{n}m_{n}<\infty $, then we call $G_{m}$ a bounded Vilenkin
group. If the generating sequence $m$ is not bounded then $G_{m}$ is said to
be an unbounded Vilenkin group. \textbf{In this paper we discuss bounded
Vilenkin groups only.}

The elements of $G_{m}$ represented by sequences
\begin{equation*}
x:=(x_{0},x_{1,...,}x_{j,...}),\qquad \left( \text{ }x_{k}\in
Z_{m_{k}}\right) .
\end{equation*}

It is easy to give a base for the neighborhood of $G_{m}:$
\begin{equation*}
I_{0}\left( x\right) :=G_{m},
\end{equation*}
\begin{equation*}
I_{n}(x):=\{y\in G_{m}\mid y_{0}=x_{0},...y_{n-1}=x_{n-1}\},\text{ }(x\in
G_{m},\text{ }n\in P).
\end{equation*}
Denote $I_{n}:=I_{n}\left( 0\right) ,$ for $n\in P$ and $\overset{-}{I_{n}}%
:=G_{m}$ $\backslash $ $I_{n}$ .

Let

\begin{equation*}
e_{n}:=\left( 0,...,0,x_{n}=1,0,...\right) \in G_{m},\qquad \left( n\in
P\right) .
\end{equation*}

Denote
\begin{equation*}
I_{N}^{k,\text{ }l}:=\left\{
\begin{array}{l}
\text{ }I_{N}(0,...,0,x_{k}\neq 0,0,...,0,x_{l}\neq 0,x_{l+1,...,\text{ }%
}x_{N-1\text{ }}),\text{ }k<l<N, \\
\text{ }I_{N}(0,...,0,x_{k}\neq 0,0,...,,0),\text{ }l=N.
\end{array}
\text{ }\right.
\end{equation*}
and
\begin{equation*}
I_{N}^{k,\alpha ,l,\beta }:=I_{N}(0,...,0,x_{k}=\alpha ,0,...,0,x_{l}=\beta
,x_{l+1,...,\text{ }}x_{N-1\text{ }}),\text{ }k<l<N.
\end{equation*}

It is evident
\begin{equation}
I_{N}^{k,\text{ }l}=\overset{m_{k}-1}{\underset{\alpha =1}{\cup }}\overset{%
m_{l}-1}{\underset{\beta =1}{\cup }}I_{N}^{k,\alpha ,l,\beta }  \label{1}
\end{equation}
and
\begin{equation}
\overset{-}{I_{N}}=\left( \overset{N-2}{\underset{k=0}{\bigcup }}\overset{N-1%
}{\underset{l=k+1}{\bigcup }}I_{N}^{k,\text{ }l}\right) \cup \left(
\underset{k=0}{\bigcup\limits^{N-1}}I_{N}^{k,\text{ }N}\right) .  \label{2}
\end{equation}

If we define the so-called generalized number system based on $m$ in the
following way :
\begin{equation*}
M_{0}:=1,\text{ \qquad }M_{k+1}:=m_{k}M_{k\text{ }},\ \qquad (k\in P).
\end{equation*}
then every $n\in P$ can be uniquely expressed as $n=\overset{\infty }{%
\underset{k=0}{\sum }}n_{j}M_{j},$ where $n_{j}\in Z_{m_{j}}$ $~(j\in P)$
and only a finite number of $n_{j}`$s differ from zero. Let $\left| n\right|
:=\max $ $\{j\in P,$ $n_{j}\neq 0\}.$

Denote by $L_{1}\left( G_{m}\right) $ the usual (one dimensional) Lebesque
space.

Next, we introduce on $G_{m}$ an ortonormal system which is called the
Vilenkin system.

At first define the complex valued function $r_{k}\left( x\right)
:G_{m}\rightarrow C,$ the generalized Rademacher functions as
\begin{equation*}
r_{k}\left( x\right) :=\exp \left( 2\pi ix_{k}/m_{k}\right) ,\text{ \qquad }%
\left( i^{2}=-1,\text{ }x\in G_{m},\text{ }k\in P\right) .
\end{equation*}

Now define the Vilenkin system $\psi :=(\psi _{n}:n\in P)$ on $G_{m}$ as:
\begin{equation*}
\psi _{n}(x):=\overset{\infty }{\underset{k=0}{\Pi }}r_{k}^{n_{k}}\left(
x\right) ,\text{ \qquad }\left( n\in P\right) .
\end{equation*}

Specifically, we call this system the Walsh-Paley one if m=2.

The Vilenkin system is ortonormal and complete in $L_{2}\left( G_{m}\right)
\,$\cite{AVD,Vi}.

Now we introduce analogues of the usual definitions in Fourier-analysis.

If $f\in L_{1}\left( G_{m}\right) $ we can establish the the Fourier
coefficients, the partial sums of the Fourier series, the Fejér means, the
Dirichlet and Fejér kernels with respect to the Vilenkin system $\psi $ in
the usual manner:
\begin{eqnarray*}
\widehat{f}\left( k\right) &:&=\int_{G_{m}}f\overline{\psi }_{k}d\mu ,\,%
\text{\qquad }\left( \text{ }k\in P\text{ }\right) , \\
S_{n}f &:&=\sum_{k=0}^{n-1}\widehat{f}\left( k\right) \psi _{k},\text{
\qquad }\left( \text{ }n\in P_{+},\text{ }S_{0}f:=0\right) , \\
\sigma _{n}f &:&=\frac{1}{n}\sum_{k=0}^{n-1}S_{k}f,\text{ \qquad }\left(
\text{ }n\in P_{+}\text{ }\right) , \\
D_{n} &:&=\sum_{k=0}^{n-1}\psi _{k\text{ }},\text{ \qquad }\left( \text{ }%
n\in P_{+}\text{ }\right) , \\
K_{n} &:&=\frac{1}{n}\overset{n-1}{\underset{k=0}{\sum }}D_{k},\text{%
\thinspace \qquad }\left( \text{ }n\in P_{+}\text{ }\right) ,
\end{eqnarray*}

Recall that
\begin{equation}
\quad \hspace*{0in}D_{M_{n}}\left( x\right) =\left\{
\begin{array}{l}
\text{ }M_{n},\text{\thinspace \thinspace \thinspace \thinspace if\thinspace
\thinspace }x\in I_{n}, \\
\text{ }0,\text{\thinspace \thinspace \thinspace \thinspace \thinspace if
\thinspace \thinspace }x\notin I_{n}.
\end{array}
\right.  \label{3}
\end{equation}

It is well-known that
\begin{equation}
\sup_{n}\int_{G_{m}}\left| K_{n}\left( x\right) \right| d\mu \left( x\right)
\leq c<\infty ,  \label{4}
\end{equation}
\begin{equation}
n\left| K_{n}\left( x\right) \right| \leq c\sum_{A=0}^{\left| n\right|
}M_{A}\left| K_{M_{A}}\left( x\right) \right| .  \label{5}
\end{equation}
\vspace{0pt}

The norm (or quasinorm) of the space $L_{p}(G_{m})$ is defined by \qquad
\qquad \thinspace \
\begin{equation*}
\left\| f\right\| _{p}:=\left( \int_{G_{m}}\left| f(x)\right| ^{p}d\mu
(x)\right) ^{1/p},\qquad \left( 0<p<\infty \right) .
\end{equation*}

The $\sigma -$algebra generated by the intervals $\left\{ I_{n}\left(
x\right) :x\in G_{m}\right\} $ will be denoted by $\digamma _{n}$ $\left(
n\in P\right) .$ Denote by $f=\left( f^{\left( n\right) },n\in P\right) $ a
martingale with respect to $\digamma _{n}$ $\left( n\in P\right) .$ (for
details see e.g. \cite{We1}). The maximal function of a martingale $f$ is
defend by \qquad
\begin{equation*}
f^{*}=\sup_{n\in P}\left| f^{\left( n\right) }\right| .
\end{equation*}

In case $f\in L_{1},$ the maximal functions are also be given by
\begin{equation*}
f^{*}\left( x\right) =\sup_{n\in P}\frac{1}{\left| I_{n}\left( x\right)
\right| }\left| \int_{I_{n}\left( x\right) }f\left( u\right) \mu \left(
u\right) \right| .
\end{equation*}

For $0<p<\infty ,$ the Hardy martingale spaces $H_{p}$ $\left( G_{m}\right) $
consist of all martingales for which
\begin{equation*}
\left\| f\right\| _{H_{p}}:=\left\| f^{*}\right\| _{L_{p}}<\infty .
\end{equation*}

If $f\in L_{1},$then it is easy to show that the sequence $\left(
S_{M_{n}}\left( f\right) :n\in P\right) $ is a martingale. If $f=\left(
f^{\left( n\right) },n\in P\right) $ is martingale then the Vilenkin-Fourier
coefficients must be defined in a slightly different manner: $\qquad \qquad $
\begin{equation*}
\widehat{f}\left( i\right) :=\lim_{k\rightarrow \infty
}\int_{G_{m}}f^{\left( k\right) }\left( x\right) \overline{\Psi }_{i}\left(
x\right) d\mu \left( x\right) .
\end{equation*}
\qquad \qquad \qquad \qquad

The Vilenkin-Fourier coefficients of $f\in L_{1}\left( G_{m}\right) $ are
the same as those of the martingale $\left( S_{M_{n}}\left( f\right) :n\in
P\right) $ obtained from $f$ .

For the martingale $f$ we consider maximal operators
\begin{eqnarray*}
\sigma ^{*}f &=&\sup_{n\in P}\left| \sigma _{n}f\right| , \\
\widetilde{\sigma }^{*}f &:&=\sup_{n\in P}\frac{\left| \sigma _{n}f\right| }{%
\log ^{2}\left( n+1\right) }.\qquad
\end{eqnarray*}

A bounded measurable function $a$ is p-atom, if there exist a dyadic
interval $I$, such that \qquad
\begin{equation*}
\left\{
\begin{array}{l}
a)\qquad \int_{I}ad\mu =0, \\
b)\ \qquad \left\| a\right\| _{\infty }\leq \mu \left( I\right) ^{-1/p}, \\
c)\qquad \text{supp}\left( a\right) \subset I.\qquad
\end{array}
\right.
\end{equation*}

\section{Formulation of Main Results}

\begin{theorem}
The \bigskip maximal operator
\end{theorem}

\begin{equation*}
\overset{\sim }{\sigma }^{*}f:=\sup_{n\in P}\frac{\left| \sigma _{n}f\right|
}{\log ^{2}\left( n+1\right) }
\end{equation*}
is bounded from the Hardy space $H_{1/2}\left( G_{m}\right) $ to the space $%
L_{1/2}\left( G_{m}\right) .$

\begin{theorem}
Let $\varphi :P_{+}\rightarrow [1,$ $\infty )$ be a nondecreasing function
satisfying the condition
\end{theorem}

\begin{equation*}
\overline{\lim_{n\rightarrow \infty }}\frac{\log ^{2}\left( n+1\right) }{%
\varphi \left( n\right) }=+\infty .
\end{equation*}
Then there exists a martingale $f\in H_{1/2},$ such that
\begin{equation*}
\sup_{n\in P}\left\| \frac{\sigma _{n}f}{\varphi \left( n\right) }\right\|
_{L_{1/2}}=\infty .
\end{equation*}

\begin{corollary}
Under condition (\ref{cond}) the maximal operator
\end{corollary}

\begin{equation*}
\sup_{n\in P}\frac{\left| \sigma _{n}f\right| }{\varphi \left( n\right) }
\end{equation*}
is not bounded from the Hardy space $H_{1/2}$ to the space $L_{1/2}.$

\section{AUXILIARY PROPOSITIONS}

\begin{lemma}
\cite{We3} Suppose that an operator $T$ is sublinear and for some $0<p\leq 1$
\end{lemma}

\begin{equation*}
\int\limits_{\overset{-}{I}}\left| Ta\right| ^{p}d\mu \leq c_{p}<\infty
\end{equation*}
for every $p$-atom $a$, where $I$ denote the support of the atom. If $T$ is
bounded from $L_{\infty \text{ }}$ to $L_{\infty },$ then
\begin{equation*}
\left\| Tf\right\| _{L_{p}\left( G_{m}\right) }\leq c_{p}\left\| f\right\|
_{H_{p}\left( G_{m}\right) }.
\end{equation*}
\bigskip

\begin{lemma}
\cite{BGG,GoAMH} Let $2<A\in P_{+},$ $k\leq s<A$ and $%
q_{A}=M_{2A}+M_{2A-2}+...+M_{2}+M_{0}.$ Then
\begin{equation*}
q_{A-1}\left| K_{q_{A-1}}(x)\right| \geq \frac{M_{2k}M_{2s}}{4},
\end{equation*}
\end{lemma}

for
\begin{equation*}
x\in I_{2A}\left( 0,...,x_{2k}\neq 0,0,...,0,x_{2s}\neq
0,x_{2s+1},...x_{2A-1}\right) ,
\end{equation*}

\begin{equation*}
k=0,\text{ }1,...,\text{ }A-3,\qquad s=k+2,\text{ }k+3,...,\text{ }A-1.
\end{equation*}

\begin{lemma}
\cite{gat} Let $A>t,$ $t,A\in P,$ $z\in I_{t}\backslash $ $I_{t+1}$ . Then
\end{lemma}

$\quad \hspace*{0in}$
\begin{equation*}
K_{M_{A}}\left( z\right) =\left\{
\begin{array}{c}
\text{ }0,\text{\qquad if }z-z_{t}e_{t}\notin I_{A}, \\
\text{ }\frac{M_{t}}{1-r_{t}\left( z\right) },\text{\qquad if }%
z-z_{t}e_{t}\in I_{A}.
\end{array}
\right.
\end{equation*}

\begin{lemma}
Let $x\in I_{N}^{k,l}$ $,$\qquad $k=0,...,N-1,$ $l=k+1,...,N.$ Then
\end{lemma}

\begin{equation*}
\int_{I_{N}}\left| K_{n}\left( x-t\right) \right| d\mu \left( t\right) \leq
\frac{cM_{l}M_{k}}{M_{N}^{2}},\,\,\,\,\text{when\thinspace \thinspace
\thinspace }n\geq M_{N}.
\end{equation*}
\textbf{Proof.} Let $x\in I_{N}^{k,\alpha ,l,\beta }$. Then applying lemma 3
we have
\begin{equation*}
K_{M_{A}}\left( x\right) =0,\,\,\text{when \thinspace \thinspace }A>l.
\end{equation*}
Hence we can suppose that $A\leq l$.

Let $k<A\leq l$. Then we have
\begin{equation}
\left| K_{M_{A}}\left( x\right) \right| =\frac{M_{k}}{\left| 1-\text{ }%
r_{k}\left( x\right) \right| }\leq \frac{m_{k}M_{k}}{2\pi \text{ }\alpha }.
\label{6}
\end{equation}

Let $x\in I_{N}^{k,\text{ }l},$ $0\leq k<l\leq N-1$ and $t\in I_{N}.$ Since $%
x-t\in $ $I_{N}^{k,\text{ }l}$ and $n\geq M_{N}$ combining (\ref{1}) (\ref{5}%
) and (\ref{6}) we obtain

\begin{equation}
n\left| K_{n}\left( x\right) \right| \leq c\overset{l-1}{\underset{A=0}{\sum
}}M_{A}M_{k}\leq cM_{k}M_{l}  \label{7}
\end{equation}

and
\begin{equation}
\int_{I_{N}}\left| K_{n}\left( x-t\right) \right| d\mu \left( t\right) \leq
\frac{cM_{k}M_{l}}{M_{N}^{2}}.  \label{8}
\end{equation}

Let $x\in I_{N}^{k,N}$ , then applying (\ref{5}) we have

\begin{equation}
\int_{I_{N}}n\left| K_{n}\left( x-t\right) \right| d\mu \left( t\right) \leq
\underset{A=0}{\overset{\left| n\right| }{\sum }}M_{A}\int_{I_{N}}\left|
K_{M_{A}}\left( x-t\right) \right| d\mu \left( t\right) .  \label{9}
\end{equation}

Let

\begin{equation*}
\left\{
\begin{array}{l}
x=\left( 0,...,0,x_{k}\neq 0,0,...0,x_{N},x_{N+1},x_{q},...,x_{\left|
n\right| -1},,...\right) , \\
t=\left( 0,...,0,x_{N},...x_{q-1},t_{q}\neq x_{q},t_{q+1},...,t_{\left|
n\right| -1},...\right) ,\,\,q=N,...,\left| n\right| -1.\text{ }
\end{array}
\right.
\end{equation*}

Using Lemma 3 in (\ref{9}) it is easy to show that

\begin{eqnarray}
&&\int_{I_{N}}\left| K_{n}\left( x-t\right) \right| d\mu \left( t\right)
\label{10} \\
&\leq &\frac{c}{n}\underset{A=0}{\overset{q-1}{\sum }}M_{A}%
\int_{I_{N}}M_{k}d\mu \left( t\right)  \notag \\
&\leq &\frac{cM_{k}M_{q}}{nM_{N}}\leq \frac{cM_{k}}{M_{N}}.  \notag
\end{eqnarray}

Let

\begin{equation*}
\left\{
\begin{array}{l}
\text{ }x=\left( 0,...,0,x_{m}\neq
0,0,...,0,x_{N},x_{N+1},x_{q},...,x_{\left| n\right| -1},...\right) ,\text{ }
\\
t=\left( 0,0,...,x_{N},...,x_{_{\left| n\right| -1}},t_{|n|},...\right) .
\end{array}
\right. .
\end{equation*}

If we apply Lemma 3 in (\ref{9}) we obtain
\begin{eqnarray}
&&\int_{I_{N}}\left| K_{n}\left( x-t\right) \right| d\mu \left( t\right)
\label{11} \\
&\leq &\frac{c}{n}\overset{\left| n\right| -1}{\underset{A=0}{\sum }}%
M_{A}\int_{I_{N}}M_{k}d\mu \left( t\right) \leq \frac{cM_{k}}{M_{N}}.  \notag
\end{eqnarray}

Combining (\ref{8}), (\ref{10}) and (\ref{11}) we complete the proof of
lemma 4.

\section{Proof of the Theorems}

\textbf{Proof of Theorem 1. }By Lemma 1, the proof of theorem 1 will be
complete, if we show that

\begin{equation*}
\int\limits_{\overline{I}_{N}}\left( \underset{n\in P}{\sup }\frac{\left|
\sigma _{n}a\right| }{\log ^{2}\left( n+1\right) }\right) ^{1/2}d\mu \leq
c<\infty ,
\end{equation*}
for every 1/2-atom $a,$ where $I$ denotes the support of the atom and
bounded from $L_{\infty }$ to $L_{\infty }.$ The boundedness follows from (%
\ref{4}).

Let $a$ be an arbitrary 1/2-atom with support$\ I$ and $\mu \left( I\right)
=M_{N}^{-1}.$ We may assume that $I=I_{N}$ $.$ It is easy to see that $%
\sigma _{n}\left( a\right) =0$ when $n\leq M_{N}$. Therefore we can suppose
that $n>M_{N}$.

Since $\left\| a\right\| _{\infty }\leq cM_{N}^{2}$ we can write
\begin{eqnarray*}
&&\frac{\left| \sigma _{n}\left( a\right) \right| }{\log ^{2}\left(
n+1\right) } \\
&\leq &\frac{1}{\log ^{2}\left( n+1\right) }\int_{I_{N}}\left| a\left(
t\right) \right| \left| K_{n}\left( x-t\right) \right| d\mu \left( t\right)
\\
&\leq &\frac{\left\| a\right\| _{\infty }}{\log ^{2}\left( n+1\right) }%
\int_{I_{N}}\left| K_{n}\left( x-t\right) \right| d\mu \left( t\right) \\
&\leq &\frac{cM_{N}^{2}}{\log ^{2}\left( n+1\right) }\int_{I_{N}}\left|
K_{n}\left( x-t\right) \right| d\mu \left( t\right) .
\end{eqnarray*}

Let $x\in I_{N}^{k,\text{ }l},\,0\leq k<l\leq N.$ Then from Lemma 4 we get
\begin{equation}
\frac{\left| \sigma _{n}\left( a\right) \right| }{\log ^{2}\left( n+1\right)
}\leq \frac{cM_{N}^{2}}{N^{2}}\frac{M_{l}M_{k}}{M_{N}^{2}}=\frac{cM_{l}M_{k}%
}{N^{2}}.  \label{12}
\end{equation}
Combining (\ref{2}) and (\ref{12}) we obtain
\begin{eqnarray*}
&&\int_{\overline{I_{N}}}\left| \widetilde{\sigma }^{*}a\left( x\right)
\right| ^{1/2}d\mu \left( x\right) \\
&=&\overset{N-2}{\underset{k=0}{\sum }}\overset{N-1}{\underset{l=k+1}{\sum }}%
\sum\limits_{x_{j}=0,j\in \{l+1,...,N-1\}}^{m_{j}-1}\int_{I_{N}^{k,\text{ }%
l}}\left| \overset{\sim }{\sigma }^{*}a\left( x\right) \right| ^{1/2}d\mu
\left( x\right) \\
&&+\overset{N-1}{\underset{k=0}{\sum }}\int_{I_{N}^{k,N}}\left| \widetilde{%
\sigma }^{*}a\left( x\right) \right| ^{1/2}d\mu \left( x\right) \\
&\leq &c\overset{N-2}{\underset{k=0}{\sum }}\overset{N-1}{\underset{l=k+1}{%
\sum }}\frac{m_{l+1}...m_{N-1}}{M_{N}}\frac{\sqrt{M_{l}M_{k}}}{N} \\
&&+c\overset{N-1}{\underset{k=0}{\sum }}\frac{1}{M_{N}}\frac{\sqrt{M_{N}M_{k}%
}}{N} \\
&\leq &c\overset{N-2}{\underset{k=0}{\sum }}\overset{N-1}{\underset{l=k+1}{%
\sum }}\frac{\sqrt{M_{k}}}{N\sqrt{M_{l}}}+c\overset{N-1}{\underset{k=0}{\sum
}}\frac{1}{\sqrt{M_{N}}}\frac{\sqrt{M_{k}}}{N}\leq c<\infty .
\end{eqnarray*}
Which complete the proof of Theorem 1.

\textbf{Proof of Theorem 2.} Let$\ \left\{ \lambda _{k};k\in P_{+}\right\} $
be an increasing sequence of the positive integers such that
\begin{equation*}
\lim_{k\rightarrow \infty }\frac{\log ^{2}\left( \lambda _{k}\right) }{%
\varphi \left( \lambda _{k}\right) }=\infty .
\end{equation*}
It is evident that for every $\lambda _{k}$ there exists a positive integers
$m_{k}^{,}$ such that $q_{m_{k}^{^{\prime }}}\leq \lambda
_{k}<q_{m_{k}^{\prime }+1}<c$ $q_{m_{k}^{^{\prime }}}.$ Since $\varphi
\left( n\right) $ is nondecreasing function.We have
\begin{equation}
\overline{\underset{k\rightarrow \infty }{\lim }}\frac{\left(
m_{k}^{^{\prime }}\right) ^{2}}{\varphi \left( q_{m_{k}^{,}}\right) }\geq
c\lim_{k\rightarrow \infty }\frac{\log ^{2}\left( \lambda _{k}\right) }{%
\varphi \left( \lambda _{k}\right) }=\infty .  \label{13}
\end{equation}

Let$\ \left\{ n_{k};k\in P_{+}\right\} \subset \left\{ m_{k}^{\prime };k\in
P_{+}\right\} $ such that
\begin{equation*}
\lim_{k\rightarrow \infty }\frac{n_{k}^{2}}{\varphi \left( q_{n_{k}}\right) }%
=\infty
\end{equation*}
and

\begin{equation*}
f_{n_{k}}\left( x\right) =D_{M_{2n_{k}+1}}\left( x\right)
-D_{M_{_{2n_{k}}}}\left( x\right) .
\end{equation*}

It is evident
\begin{equation*}
\widehat{f}_{n_{k}}\left( i\right) =\left\{
\begin{array}{l}
\text{ }1,\text{ if }i=M_{_{2n_{k}}},...,M_{2n_{k}+1}-1, \\
\text{ }0,\text{otherwise.}
\end{array}
\right.
\end{equation*}
Then we can write
\begin{equation}
S_{i}\left( f_{n_{k}}(x)\right) =\left\{
\begin{array}{l}
D_{i}\left( x\right) -D_{M_{_{2n_{k}}}}\left( x\right) \text{ },\text{ if }%
i=M_{_{2n_{k}}},...,M_{2n_{k}+1}-1, \\
\text{ }f_{n_{k}}\left( x\right) \text{ },\text{if }i\geq M_{2n_{k}+1}, \\
0,\text{ \qquad otherwise}.
\end{array}
\right.  \label{14}
\end{equation}

From (\ref{3}) we get
\begin{eqnarray}
&&\left\| f_{n_{k}}\right\| _{H_{1/\text{ }2}}  \label{15} \\
&=&\left\| \sup\limits_{n\in P}S_{M_{n}}f_{n_{k}}\right\| _{L_{1/2}}  \notag
\\
&=&\left\| D_{M_{2n_{k}+1}}-D_{M_{_{2n_{k}}}}\right\| _{L_{1/2}}  \notag \\
&=&\left( \int_{I_{_{2n_{k}}}\backslash \text{ }%
I_{_{2n_{k}+1}}}M_{_{2n_{k}}}^{1/2}d\mu \left( x\right)
+\int_{I_{_{2n_{k}+1}}}\left( M_{_{2n_{k}+1}}-M_{_{2n_{k}}}\right)
^{1/2}d\mu \left( x\right) \right) ^{2}  \notag \\
&=&\left( \frac{m_{_{2n_{k}}}-1}{M_{2n_{k}+1}}M_{_{2n_{k}}}^{1/2}+\frac{%
\left( m_{_{2n_{k}}}-1\right) ^{1/2}}{M_{_{2n_{k}}+1}}M_{_{2n_{k}}}^{1/2}%
\right) ^{2}  \notag \\
&\leq &\frac{c}{M_{_{2n_{k}}}}.  \notag
\end{eqnarray}
By (\ref{14}) we can write:
\begin{eqnarray*}
&&\frac{\left| \sigma _{q_{n_{k}}}\left( f_{n_{k}}(x)\right) \right| }{%
\varphi \left( q_{n_{k}}\right) } \\
&=&\frac{1}{\varphi \left( q_{n_{k}}\right) q_{n_{k}}}\left| \overset{%
q_{n_{k}}-1}{\underset{j=0}{\sum }}S_{j}f_{n_{k}}(x)\right| \\
&=&\frac{1}{\varphi \left( q_{n_{k}}\right) q_{n_{k}}}\left| \overset{%
q_{n_{k}}-1}{\underset{j=M_{_{2n_{k}}}}{\sum }}S_{j}f_{n_{k}}(x)\right| \\
&=&\frac{1}{\varphi \left( q_{n_{k}}\right) q_{n_{k}}}\left| \overset{%
q_{n_{k}}-1}{\underset{j=M_{_{2n_{k}}}}{\sum }}\left( D_{j}\left( x\right)
-D_{M_{_{2n_{k}}}}\left( x\right) \right) \right| \\
&=&\frac{1}{\varphi \left( q_{n_{k}}\right) q_{n_{k}}}\left| \overset{%
q_{n_{k}-1}-1}{\underset{j=0}{\sum }}\left( D_{j+M_{_{2n_{k}}}}\left(
x\right) -D_{M_{_{2n_{k}}}}\left( x\right) \right) \right|
\end{eqnarray*}
Since

\begin{equation*}
D_{_{j+M_{_{2n_{k}}}}}\left( x\right) -D_{M_{_{2n_{k}}}}\left( x\right)
=\psi _{M_{_{2n_{k}}}}D_{j},\text{ }\,j=1,2,..,M_{_{2n_{k}}}-1,
\end{equation*}
we obtain
\begin{eqnarray*}
&&\frac{\left| \sigma _{q_{n_{k}}}\left( f_{n_{k}}(x)\right) \right| }{%
\varphi \left( q_{n_{k}}\right) } \\
&=&\frac{1}{\varphi \left( q_{n_{k}}\right) q_{n_{k}}}\left| \overset{%
q_{n_{k}-1}-1}{\underset{j=0}{\sum }}D_{j}\left( x\right) \right| \\
&=&\frac{1}{\varphi \left( q_{n_{k}}\right) }\frac{q_{n_{k}-1}}{q_{n_{k}}}%
\left| K_{q_{n_{k}-1}}\left( x\right) \right| .
\end{eqnarray*}

Let $x\in $ $I_{_{2n_{k}}}^{2s,2l}$. Then from Lemma 2 we obtain
\begin{equation*}
\frac{\left| \sigma _{q_{n_{k}}}\left( f_{n_{k}}(x)\right) \right| }{\varphi
\left( q_{n_{k}}\right) }\geq \frac{cM_{2s}M_{2l}}{M_{_{2n_{k}}}\varphi
\left( q_{n_{k}}\right) }.
\end{equation*}
Hence we can write:
\begin{eqnarray*}
&&\int_{G_{m}}\left| \frac{\sigma _{q_{n_{k}}}\left( f_{n_{k}}(x)\right) }{%
\varphi \left( q_{n_{k}}\right) }\right| ^{1/2}d\mu \left( x\right) \\
&\geq &\text{ }\overset{n_{k}-3}{\underset{s=0}{\sum }}\overset{n_{k}-1}{%
\underset{l=s+1}{\sum }}\overset{m_{2l+1}}{\underset{x_{2l+1=0}}{\sum }}...%
\overset{m_{2n_{k}-1}}{\underset{x_{2n_{k}-1}=0}{\sum }}%
\int_{I_{_{2n_{k}}}^{2s,2l}}\left| \frac{\sigma _{q_{n_{k}}}\left(
f_{n_{k}}(x)\right) }{\varphi \left( q_{n_{k}}\right) }\right| ^{1/2}d\mu
\left( x\right) \\
&\geq &c\overset{n_{k}-3}{\underset{s=0}{\sum }}\overset{n_{k}-1}{\underset{%
l=s+1}{\sum }}\frac{m_{_{2l+1}}...m_{2n_{k}-1}}{M_{2n_{k}}}\frac{\sqrt{%
M_{2s}M_{2l}}}{\sqrt{\varphi \left( q_{n_{k}}\right) M_{_{2n_{k}}}}} \\
&\geq &c\overset{n_{k}-3}{\underset{s=0}{\sum }}\overset{n_{k}-1}{\underset{%
l=s+1}{\sum }}\frac{\sqrt{M_{2s}}}{\sqrt{M_{2l}M_{_{2n_{k}}}\varphi \left(
q_{n_{k}}\right) }} \\
&\geq &\frac{cn_{k}}{\sqrt{M_{_{2n_{k}}}\varphi \left( q_{n_{k}}\right) }}.
\end{eqnarray*}
From (\ref{15}) we have
\begin{eqnarray*}
&&\frac{\left( \int_{G_{m}}\left| \frac{\sigma _{q_{n_{k}}}\left(
f_{n_{k}}(x)\right) }{\varphi \left( n_{k}\right) }\right| ^{1/2}d\mu \left(
x\right) \right) ^{2}}{\left\| f_{n_{k}}\right\| _{H_{1/\text{ }2}}} \\
&\geq &\frac{cn_{k}^{2}}{M_{_{2n_{k}}}\varphi \left( q_{n_{k}}\right) }%
M_{2n_{k}} \\
&\geq &\frac{cn_{k}^{2}}{\varphi \left( q_{n_{k}}\right) }\rightarrow \infty
,\qquad \text{when \qquad }k\rightarrow \infty .
\end{eqnarray*}

Theorem 2 is proved.

\end{document}